\documentclass[11pt]{article}
\usepackage{latexsym}
\usepackage{amsfonts}
\usepackage{enumerate}
\usepackage{multicol}
\usepackage{graphicx}
\usepackage{amssymb}
\usepackage{amsmath}
\usepackage{epic}

\title{More on the $h$-critical numbers of finite abelian groups}

\author{B\'{e}la Bajnok \\[.1in] Department of Mathematics, Gettysburg College \\
Gettysburg, PA 17325-1486 USA \\E-mail:  bbajnok@gettysburg.edu \\[.4in]}

\date{November 21, 2016}

\newtheorem{thm}{Theorem}

\newtheorem{prop}[thm]{Proposition}

\begin{document}

\maketitle

\begin{abstract}

For a finite abelian group $G$, a nonempty subset $A$ of $G$, and a positive integer $h$, we let $hA$ denote the $h$-fold sumset of $A$; that is, $hA$ is the collection of sums of $h$ not-necessarily-distinct elements of $A$.  Furthermore, for a positive integer $s$, we set $[0,s] A=\cup_{h=0}^s h A$.  We say that $A$ is a generating set of $G$ if there is a positive integer $s$ for which $[0,s] A=G$.

The $h$-critical number $\chi (G,h)$ of $G$ is defined as the smallest positive integer $m$ for which $hA=G$ holds for every $m$-subset $A$ of $G$; similarly, $\chi (G,[0,s])$ is the smallest positive integer $m$ for which $[0,s]A=G$ holds for every $m$-subset $A$ of $G$.  We define $\widehat{\chi}  (G, h)$  as the smallest positive integer $m$ for which $hA=G$ holds for every generating $m$-subset $A$ of $G$; $\widehat{\chi}  (G, [0,s])$ is defined similarly.

The value of $\chi (G,h)$ has been determined by this author for all $G$ and $h$, and $\widehat{\chi}  (G, [0,s])$ was introduced and resolved for some special cases by Klopsch and Lev.  Here we determine the remaining two quantities in all cases.

\end{abstract}

\noindent AMS Mathematics Subject Classification:  \\ Primary: 11B75; \\ Secondary: 05D99, 11B25, 11P70, 20K01.

\noindent Key words and phrases: \\ critical number, abelian groups, sumsets, generating sets.

\thispagestyle{empty}

\section{Introduction}

Let $G$ be a finite abelian group of order $n \geq 2$, written in additive notation.  When $G$ is cyclic, we will identify it with $\mathbb{Z}_n=\mathbb{Z}/n\mathbb{Z}$.  More generally, we recall that $G$ has a unique {\em type} $(n_1,\dots,n_r)$, where $r$ and $n_1, \dots, n_r$ are positive integers  so that $n_1 \geq 2$, $n_i$ is a divisor of $n_{i+1}$ for
 $i=1,\dots,r-1$, and $$G \cong \mathbb{Z}_{n_1} \times \cdots \times \mathbb{Z}_{n_r};$$ here $r$ is the {\em rank} of $G$ and $n_r$ is the {\em exponent} of $G$.  

For a nonempty subset $A$ of $G$, we let $\Sigma A$ be the set of all subset sums of $A$; that is, $$\Sigma A = \left \{ \Sigma_{a \in B} a \mid B \subseteq A \right \}$$  (with the subset sum of the empty-set defined as $0$).  We then let $\mathrm{cr}(G)$ denote the  smallest integer $m$ for which $\Sigma A = G$ holds for all $m$-subsets $A$ of $G$; the smallest integer $m$ for which $\Sigma A = G$ holds for all $m$-subsets $A$ of $G \setminus \{0\}$ is denoted by $\mathrm{cr}^*(G)$.

The study of {\em critical numbers} originated with the 1964 paper \cite{ErdHei:1964a} of Erd\H{o}s and Heilbronn, in which they asked for $\mathrm{cr}^*(\mathbb{Z}_p)$ for prime values of $p$.  
It took nearly half a century, but now, due to the combined results of Diderrich and Mann \cite{DidMan:1973a}, Diderrich \cite{Did:1975a},  Mann and Wou \cite{ManWou:1986a}, Dias Da Silva and Hamidoune \cite{DiaHam:1994a}, Gao and Hamidoune \cite{GaoHam:1999a}, Griggs \cite{Gri:2001a}, and  Freeze, Gao, and Geroldinger \cite{FreGaoGer:2009a, FreGaoGer:2014a}, we have the critical number of every group:

\begin{thm} [The combined results of authors above] \label{combined critical}

Suppose that $G$ is an abelian group of order $n \geq 10$, and let $p$ be the smallest prime divisor of $n$.  Then
$$\mathrm{cr}^*(G)=\mathrm{cr}(G)-1=\left\{
\begin{array}{ll}
\lfloor 2 \sqrt{n-2} \rfloor & \mbox{if $G$ is cyclic of order $n=p$ or $n=pq$ where} \\
& \mbox{$q$ is prime and $3 \leq p \leq q \leq p+\lfloor 2 \sqrt{p-2} \rfloor+1$}\footnotemark, \\ \\ 
n/p+p-2 & \mbox{otherwise}.
\end{array}
\right.
$$

\end{thm}
\footnotetext{Observe that $\lfloor 2 \sqrt{n-2} \rfloor=n/p+p-1$ in this case.}
We note that, while it is easy to see that $\mathrm{cr}(G)$ is at least one more than $\mathrm{cr}^*(G)$, there is no obvious reason known for the fact that they differ by exactly one. 

In this paper we consider some variations of the critical numbers defined above.  

We recall the following definitions.  For a positive integer $h$ and a nonempty subset $A$ of $G$, we let $hA$ denote the {\em $h$-fold sumset} of $A$; that is, $hA$ is the collection of sums of $h$ not-necessarily-distinct elements of $A$.  Additionally, for a positive integer $s$, we set $[0,s] A=\cup_{h=0}^s h A$.  Recall also that, for a subset $A$ of $G$, $\langle A \rangle$ is the subgroup generated by $A$ in $G$; that is, $\langle A \rangle$ is the intersection of all subgroups $H$ of $G$ for which $A \subseteq H$.  When $\langle A \rangle=G$, we say that $A$ is a {\em generating set} of $G$.  Clearly, $A$ is a generating set of $G$ if, and only if, there is a positive integer $s$ for which $[0,s] A=G$.

The subject of our paper is the study of the following four quantities:
\begin{eqnarray*}
{\chi}  (G, h) & = & \min \{m \; : \; A \subseteq G,  |A| \geq m \Rightarrow h  A=G \},\\
{\chi}  (G, [0,s]) & = &  \min \{m \; : \; A \subseteq G,  |A| \geq m \Rightarrow [0,s]  A=G \},\\
\widehat{\chi}  (G, h) & = &  \min \{m \; : \; A \subseteq G, \langle A \rangle =G,  |A| \geq m \Rightarrow h  A=G \},\\
\widehat{\chi}  (G, [0,s]) & = &  \min \{m \; : \; A \subseteq G, \langle A \rangle =G,  |A| \geq m \Rightarrow [0,s]  A=G \}.
\end{eqnarray*}
It is easy to see that for all $G$ and $h$ we have $hG=G$, so all four quantities are well defined.

The value of  ${\chi}  (G, h)$ is now known for every $G$ and $h$.  To state the result, we let $D(n)$ denote the set of positive divisors of the order $n$ of $G$; then  set
$$f_d(n,h)=\left( \left \lfloor \frac{d-2}{h} \right \rfloor +1  \right) \cdot \frac{n}{d}$$
for each $d \in D(n)$, and let
$$v(n,h)= \max \left\{ f_d(n,h)  \; : \; d \in D(n) \right\}.$$ We should note that the function $v(n,h)$ has appeared elsewhere in additive combinatorics already.  For example, according to the classical result of 
Diananda and Yap (see \cite{DiaYap:1969a}), the maximum size of a sum-free set (that is, a set $A$ that is disjoint from $2A$) in the cyclic group $\mathbb{Z}_n$ is given by
$$v(n,3)= \left\{
\begin{array}{ll}
\left(1+\frac{1}{p}\right) \frac{n}{3} & \mbox{if $n$ has prime divisors congruent to 2 mod 3,} \\ & \mbox{and $p$ is the smallest such divisor,}\\ \\
\left\lfloor \frac{n}{3} \right\rfloor & \mbox{otherwise;}\\
\end{array}\right.$$
see e.g.~\cite{Baj:2009a} for generalizations.  

Our result for ${\chi}  (G, h)$ is then the following:

\begin{thm} [Bajnok; cf.~\cite{Baj:2014a}]  \label{thm bajnok chi h}
For any abelian group $G$ of order $n$ and for all positive integers $h$ we have $${\chi}  (G, h)=v(n,h)+1.$$

\end{thm}    
By Theorem \ref{thm bajnok chi h}, the size of the largest {\em $h$-incomplete subset} of any group of order $n$---that is, a subset whose $h$-fold sumset is not the entire group---equals $v(n,h)$.

Given Theorem \ref{thm bajnok chi h}, the evaluation of ${\chi}  (G, [0,s])$ is immediate.  Clearly, ${\chi}  (G, s)$ is an upper bound for 
 ${\chi}  (G, [0,s])$.  Suppose then that $A$  is an $s$-incomplete subset of $G$ of size ${\chi}  (G, s)-1$.  Choose an element $a_0 \in A$, and let
$$B=A-a_0=\{a-a_0 \mid a \in A\}.$$  Since $|A|=|B|$ and $|sA|=|sB|$, $B$ is also an $s$-incomplete subset of $G$ of size ${\chi}  (G, s)-1$.  But $[0,s]B=sB$, and thus $[0,s]B \neq G$, which implies that 
${\chi}  (G, [0,s])$ is an upper bound for 
 ${\chi}  (G, s)$.  Therefore:

\begin{thm}   \label{thm bajnok chi [0,s]}
For any abelian group $G$ of order $n$ and for all positive integers $s$ we have $${\chi}  (G, [0,s])={\chi} (G,s)=v(n,s)+1.$$

\end{thm}

The quantity $\widehat{\chi}  (G, [0,s])$ was introduced and investigated by Klopsch and Lev in \cite{KloLev:2009a} (though earlier works had treated the case of elementary abelian 2-groups).  The value of $\widehat{\chi}  (G, [0,s])$ is not known in general.   In Section \ref{chi cap G [0,s]} we provide a general lower bound for $\widehat{\chi}  (G, [0,s])$, and summarize the main results in a way that allows for comparisons to our other quantities and may generate renewed interest.

Finally, we consider $\widehat{\chi}  (G, h)$.  In Section \ref{chi cap G h} we prove that, perhaps surprisingly, ${\chi}  (G, h)$ and $\widehat{\chi}  (G, h)$ are always equal:

\begin{thm} \label{thm bajnok chi cap  h}
For any abelian group $G$ of order $n$ and for all positive integers $h$ we have $$\widehat{\chi} (G, h)={\chi}  (G, h)=v(n,h)+1.$$

\end{thm}

\section{On the value of $\widehat{\chi}  (G, [0,s])$} \label{chi cap G [0,s]}

It is an easy exercise to show that, if $G$ is an abelian group of order $n$ that is not isomorphic to $\mathbb{Z}_2$, then $$\widehat{\chi}  (G,[0,1])={\chi}  (G,[0,1])=v(n,1)+1=n,$$ and if it is not isomorphic to $\mathbb{Z}_2$ or $\mathbb{Z}_2^2$, then $$\widehat{\chi}  (G,[0,2])={\chi}  (G,[0,2])=v(n,2)+1=\lfloor n/2 \rfloor +1.$$  For $s =3$, the result is considerably more complicated; in particular, it depends on the structure of $G$ and not just on the order $n$ of $G$:  

\begin{thm} [Klopsch and Lev; cf.~\cite{KloLev:2009a}]  \label{thm Klopsch Lev s=3 result}
If $G$ is a finite abelian group of order $n$ that is not isomorphic to an elementary abelian 2-group, then   
$$\widehat{\chi}  (G,[0,3])=\left\{
  \begin{array}{cl}
\left( 1+ \frac{1}{d} \right) \cdot \frac{n}{3}+1 & \mbox{if} \; \mbox{$G$ has a subgroup whose order is congruent to 2 mod 3}\\
& \mbox{that is not isomorphic to an elementary abelian 2-group,} \\
& \mbox{and $d$ is the minimum size of such a subgroup}; \\ \\
\left \lfloor \frac{n}{3}\right \rfloor +1& otherwise.
\end{array}
\right.$$

\end{thm} 
(The value of $\widehat{\chi}  (G,[0,s])$ for elementary abelian 2-groups had been determined earlier---see below.   
We should note that this result appeared in \cite{KloLev:2009a} via different expressions.)  Theorems \ref{thm bajnok chi [0,s]} and \ref{thm Klopsch Lev s=3 result} allow for an interesting comparison; for example, we see that 
$$\widehat{\chi}  (G,[0,3])={\chi}  (G,[0,3])$$ holds if, and only if, $n$ is odd.  

The authors of \cite{KloLev:2009a} warn that an expression for $\widehat{\chi}  (G,[0,s])$ when $s \geq 4$ ``is very difficult, if at all feasible.''  One thus may focus on special types of groups.  It appears that only two such results are currently known: 

\begin{thm} [Klopsch and Lev; cf.~\cite{KloLev:2009a}] \label{Klo Lev cyclic}
Let $n$ and $s$ be positive integers.  If $n \leq s+1$, then $\widehat{\chi}  (\mathbb{Z}_n,[0,s])=1$; otherwise, we have
$$\widehat{\chi}  (\mathbb{Z}_n,[0,s])=\max \left\{ f_d(n,s) \; : \; d \in D(n), d \geq s+2 \right\}+1,$$ where $f_d(n,s)$ is as defined above.

\end{thm}

\begin{thm} [Lev; cf.~\cite{Lev:2003a}] \label{Lev 2-group}
Let $r$ and $s$ be positive integers, $s \geq 2$.  If $r \leq s$, then $\widehat{\chi}  (\mathbb{Z}_2^r,[0,s])=1$; otherwise we have
$$\widehat{\chi}  (\mathbb{Z}_2^r,[0,s])=(s+2) \cdot 2^{r-s-1}+1.$$

\end{thm}

While the two results appear dissimilar, the following general lower bound evaluates to the stated values of $\widehat{\chi}  (G,[0,s])$ in both cases.

\begin{prop} \label{chi cap lower [0,s]}
Let $G$ be an abelian group of order $n$, and let $H$ be a subgroup of $G$ of index $d>1$ for which $G/H$ is of type $(d_1,\dots,d_t)$.  For each $i=1,\dots, t$, let $c_i$ be a positive integer with $c_i \leq d_i-1$, and suppose that
$$\Sigma_{i=1}^t \left \lceil (d_i-1)/c_i \right \rceil \geq s+1.$$  Then we have
$$\widehat{\chi}  (G,[0,s]) \geq \left ( 1+ \Sigma_{i=1}^t c_i \right) \cdot n/d +1.$$

\end{prop}

Before proving Proposition \ref{chi cap lower [0,s]}, we deduce how it provides exact lower bounds for $\widehat{\chi}  (G,[0,s])$ in both Theorems \ref{Klo Lev cyclic} and \ref{Lev 2-group}.

Suppose first that $G$ is cyclic and of order $n$.  Clearly, if $n \leq s+1$, then $\widehat{\chi}  (G,[0,s])=1$.  Assume then that $n \geq s+2$, and let $d$ be any divisor of $n$ with $d \geq s+2$.   Let $H$ be a subgroup of $G$ of index $d$, in which case $G/H$ is cyclic and of order $d$.  Then, with $t=1$ and $c=\left \lfloor (d-2)/s \right \rfloor$, we have $c \geq 1$ and $\left \lceil (d-1)/c \right \rceil \geq s+1,$ so by Proposition \ref{chi cap lower [0,s]}, we get
$$\widehat{\chi}  (G,[0,s]) \geq \left( \left \lfloor \frac{d-2}{s} \right \rfloor +1  \right) \cdot \frac{n}{d} +1=f_d(n,s)+1,$$ and our claim follows.

Next, we consider $\mathbb{Z}_2^r$, the elementary abelian 2-group of rank $r$.  The result is trivial when $r \leq s$, so assume that $s+1 \leq r$, and let $t$ be an integer with $$s+1 \leq t \leq r.$$  Then choosing $H=\mathbb{Z}_2^t$ and $c_i=1$ for all $i \in \{1,\dots,t\}$, 
Proposition \ref{chi cap lower [0,s]} implies that 
$$\widehat{\chi}  (\mathbb{Z}_2^r, [0,s]) \geq  (t+1) \cdot 2^{r-t}+1;$$ in particular,
we have
$$\widehat{\chi}  (\mathbb{Z}_2^r, [0,s]) \geq  (s+2) \cdot 2^{r-s-1}+1,$$ as claimed.

{\em Proof of Proposition \ref{chi cap lower [0,s]}.}  We shall prove our claim by exhibiting a subset $A$ of $G$ of size 
$$|A|=\left ( 1+ \Sigma_{i=1}^t c_i \right) \cdot n/d$$ for which $\langle A \rangle =G$, but $[0,s]A \neq G$.

Let us identify $G/H$ with 
$$K=\mathbb{Z}_{d_1} \times \cdots \times \mathbb{Z}_{d_t},$$ and for $i=1,\dots,t$, set $$N_i=\{1,2, \dots, c_i \} \subseteq \mathbb{Z}_{d_i}.$$  Now let 
$$B_i = \{0\}^{i-1} \times N_i \times \{0\}^{t-i}$$ with the understanding that $\{0\}^0$ is to be ignored, and let $B_0=\{0\}^t$.  

Consider $B=\cup_{i=0}^t B_i.$
We have $\langle B \rangle =K$, and  $$|B|=1+\Sigma_{i=1}^t c_i.$$
Furthermore, observe that when
$$\Sigma_{i=1}^t \left \lceil (d_i-1)/c_i \right \rceil \geq s+1,$$ then $[0,s]B=sB \neq K.$

Now set $A =\pi^{-1}(B)$, where $\pi: G \rightarrow G/H$ is the canonical homomorphism; it is easy to see that $A$ satisfies our requirements.  $\Box$

\section{The evaluation of $\widehat{\chi}  (G, h)$} \label{chi cap G h}

In this section we prove Theorem \ref{thm bajnok chi cap  h}, namely, that for any abelian group $G$ of order $n$ and for all positive integers $h$, we have $$\widehat{\chi} (G, h)={\chi}  (G, h)=v(n,h)+1.$$

Note that, obviously, 
$$\widehat{\chi} (G, h) \leq {\chi} (G, h),$$ so by Theorem \ref{thm bajnok chi h}, we have $$ \widehat{\chi} (G, h) \leq v(n,h)+1.$$
Therefore, to establish Theorem \ref{thm bajnok chi cap  h}, it suffices to find a subset $A$ of $G$ of size $v(n,h)$ for which $\langle A \rangle =G$, but  $hA \neq G$.

We proceed by induction on the order $n$ of $G$.  The claim can be easily verified for $n=2$; we will assume that it also holds for all groups of order at most $n-1$. 

Recall that we have set $$v(n,h)= \max \left\{ f_d(n,h)  \; : \; d \in D(n) \right\},$$ where $D(n)$ is the set of positive divisors of $n$, and 
$$f_d(n,h)= \left( \left \lfloor \frac{d-2}{h} \right \rfloor +1  \right) \cdot \frac{n}{d}.$$

We consider two cases.

{\em Case 1}:  There is a $d_0 \in D(n) \setminus \{n\}$ for which $v(n,h)=f_{d_0}(n,h)$. 

Let $H$ be a subgroup of index $d_0$ in $G$.  Then $G/H$ has order $d_0<n$, so by our inductive hypotheses, it contains a subset $B$ of size $v(d_0,h)$ for which $\langle B \rangle =G/H$, but $ hB \neq G/H$.  

Let $\pi: G \rightarrow G/H$ denote the canonical homomorphism, and let $A= \pi^{-1} (B)$.  Then $\langle A \rangle =G$,  $ hA \neq G$, and 
the size of $A$ is $v(n,h)$, since
$$|A|=|B| \cdot |H|=v(d_0,h) \cdot \frac{n}{d_0} \geq f_{d_0}(d_0,h) \cdot \frac{n}{d_0}=\left( \left \lfloor \frac{d_0-2}{h} \right \rfloor + 1 \right) \cdot \frac{n}{d_0} =f_{d_0}(n,h)= v(n,h).$$  This completes the proof of Case 1.

{\em  Case 2}:  For all $d \in D(n) \setminus \{n\}$, $$f_d(n,h) \leq v(n,h)-1.$$
We then must have $$v(n,h)=f_n(n,h)=\left \lfloor \frac{n-2}{h} \right \rfloor +1.$$

Let us consider first the case when $n$ equals a prime number $p$; we will then identify $G$ with the cyclic group $\mathbb{Z}_p$.    
Let $$A=\left \{1,2, \dots, \left \lfloor \frac{p-2}{h} \right \rfloor +1\right \} \subseteq \mathbb{Z}_p.$$  Clearly, $A$ has size $v(p,h)$, and $\langle A \rangle =\mathbb{Z}_p$.  Moreover,
$$hA=\left \{h,h+1, \dots, h \cdot \left \lfloor \frac{p-2}{h} \right \rfloor +h\right \},$$ from which we see that $h-1 \not \in hA$, and thus $hA \neq \mathbb{Z}_p$.  Therefore, our claim holds for prime values of $n$.

Assume now that $n$ is composite.  We will show that for each prime divisor $p$ of $n$, $p-1$ must be divisible by $h$.

To do so, let $p$ be a prime divisor of $n$, and let $$p-2=ch+r$$ for some (unique) integers $c$ and $r$ with $$0 \leq r \leq h-1.$$ Note that $n$ is composite, so $p <n$.  Therefore, by assumption, we have
$$f_p(n,h) \leq v(n,h)-1=\left \lfloor \frac{n-2}{h} \right \rfloor.$$ Here
$$f_p(n,h)=\left( \left \lfloor \frac{p-2}{h} \right \rfloor +1  \right) \cdot \frac{n}{p}=\left( \frac{p-2-r}{h} +1  \right) \cdot \frac{n}{p}=\frac{n}{h}+\frac{h-2-r}{h}\cdot \frac{n}{p},$$ which is more than 
$\left \lfloor \frac{n-2}{h} \right \rfloor$, unless $r=h-1$.  But then $$p-2=ch+h-1,$$ so $p-1$ is divisible by $h$, as claimed.

This implies that every positive divisor of $n$ is congruent to 1 mod $h$; in particular, so is $n$, and thus
$$v(n,h)=f_n(n,h)=\left \lfloor \frac{n-2}{h} \right \rfloor +1=\frac{n-1}{h}.$$
We thus need to find a subset $A$ of $G$ of size $\frac{n-1}{h}$ for which $\langle A \rangle =G$, but  $hA \neq G$.

Suppose that $G$ has rank $r$ and that it is of type $(n_1,\dots,n_r)$; that is, $$G \cong \mathbb{Z}_{n_1} \times \cdots \times \mathbb{Z}_{n_r}$$ for integers $n_1, \dots, n_r$ for which $n_1 \geq 2$ and $n_i$ is a divisor of $n_{i+1}$ for each $i=1,\dots,r-1$.  As we just proved, $n_i-1$ is divisible by $h$ for each $i=1,\dots,r$.

We construct $A$ as follows: for each $i=1,\dots,r$, let $$N_i=\left \{ 1,2,\dots,\frac{n_i-1}{h} \right\} \subseteq \mathbb{Z}_{n_i}$$ and
$$A_i=\mathbb{Z}_{n_1} \times \cdots \times \mathbb{Z}_{n_{i-1}} \times N_i \times \{0\}^{r-i};$$ we then set $$A=\cup_{i=1}^r A_i.$$
Then $$|A_i|=n_1 \cdots n_{i-1} \cdot \frac{n_i-1}{h};$$ since the $r$ sets are pairwise disjoint, the size of $A$ equals
$$|A|=\sum_{i=1}^r |A_i|= \sum_{i=1}^r n_1 \cdots n_{i-1} \cdot \frac{n_i-1}{h}=\frac{n_1\cdots n_r-1}{h}=\frac{n-1}{h}=v(n,h).$$

Furthermore, since $n_r-1$ is divisible by $h$, we have $n_r \geq h+1$, and thus 
$$A \supseteq A_r \supseteq \mathbb{Z}_{n_1} \times \cdots \times \mathbb{Z}_{n_{r-1}} \times \{1\},$$
which implies that $\langle A \rangle=G$.  

Finally, we show that $hA \neq G$ by verifying that $0 \not \in hA$.  Let $a_1, \dots, a_h$ be (not-necessarily distinct) elements of $A$.  For each $i=1,\dots,r$, there is a unique $k_i \in \{1,\dots,r\}$ for which $a_i \in A_{k_i}$; let $$k=\max \{k_i \mid i=1,\dots,r\}.$$  But then the $k$-th coordinate of $a_1+\cdots + a_h$ (as an element of $\mathbb{Z}_{n_1} \times \cdots \times \mathbb{Z}_{n_r}$) is at least 1 and at most $n_k-1,$ so $$a_1+\cdots + a_h \neq 0.$$
This completes our proof. $\Box$

\end{document}